\documentstyle[12pt,twoside]{article}
\newtheorem{theorem}{\bf Theorem}[section]
\newtheorem{proposition}[theorem]{\bf Proposition}
\newtheorem{lemma}[theorem]{\bf Lemma}
\newtheorem{corollary}[theorem]{\bf Corollary}
\newtheorem{example}[theorem]{\bf Example}

\input{amssym}

\date{}
\begin{document}

\title{{\Large\bf Weak amenability of weighted group algebras}}

\author{{\normalsize\sc M. J. Mehdipour and A. Rejali\footnote{Corresponding author}}}
\maketitle

{\footnotesize  {\bf Abstract.} In this paper,  we study weak amenability of Beurling algebras. To this end, we introduce the notion inner quasi-additive functions and prove that for a locally compact group $G$, the Banach algebra $L^1(G, \omega)$ is weakly amenable if and only if every non-inner quasi-additive function in $L^\infty(G, 1/\omega)$ is unbounded. This provides an answer to the question concerning weak amenability of $L^1(G, \omega)$ and improve some known results in connection with it.}
                     %-----------------------------------------
                     %-----------------------------------------
                     %-----------------------------------------
{}
                     %-----------------------------------------
                     %-----------------------------------------
                     %-----------------------------------------

\section{\normalsize\bf Introduction}

Throughout this paper $G$ is a locally compact group with an identity element $e$. Let $L^1(G)$ and $M(G)$ be group and measure algebras of $G$, respectively. Let also $L^\infty(G)$ be the Lebesgue space of bounded Borel measurable functions on $G$; see \cite{d,hr} for an extensive study of these spaces.

Let us recall that a continuous function $\omega: G\rightarrow [1, \infty)$ is called a \emph{weight} \emph{function} if for every $x, y\in G$
$$
\omega(xy)\leq\omega(x)\;\omega(y)\quad\hbox{and}\quad\omega(e)=1.
$$
An elementary computation shows that the functions $\omega^*$ and $\omega^\otimes $ defined by
$$
\omega^*(x):=\omega(x)\omega(x^{-1})\quad\hbox{and}\quad\omega^\otimes (x, y)=\omega(x)\omega(y)
$$
are weight functions on $G$ and $G\times G$, respectively. Let the function $\Omega: G\times G\rightarrow (0, 1]$ be defined as follows:
$$
\Omega(x, y)=\frac{\omega(xy)}{\omega^\otimes (x, y)}.
$$
Then $\Omega$ is called \emph{zero cluster} if for every pair of sequences $(x_n)_n$ and $(y_m)_m$ of distinct elements in $G$, we have
\begin{eqnarray*}\label{positive cluster}
\lim_n\lim_m\Omega(x_n, y_m)=0=\lim_m\lim_n\Omega(x_n, y_m),
\end{eqnarray*}
whenever both iterated limits exist.

Let $L^1(G, \omega)$ be the Banach space of all Borel measurable functions $f$ on $G$ such that $\omega f\in L^1(G)$. Then $L^1(G, \omega)$ with the convolution product ``$\ast$" and the norm $\|f\|_\omega=\|\omega f\|_1$ is a Banach algebras.
Let $L^\infty (G,1/\omega)$ be the space of all Borel measurable
functions $f$ on $G$ with $f/\omega\in L^\infty (G)$. Then
$L^\infty (G,1/\omega)$ with the norm $$
\|h\|_{\infty,\;\omega}=\|h/\omega\|_\infty,
$$
and the
multiplication $``\cdot_\omega"$ defined by
$$
h\cdot_\omega k=hk/\omega\quad\quad\quad (h, k\in L^\infty(G,
1/\omega))
$$
is a commutative
$C^*-$algebra. Also, $L^\infty(G, 1/\omega)$ and $L^1(G, \omega)^*$ are isometrically isomorphic with the duality given through
$$
\langle f, h\rangle=\int_G f(x) h(x)\; dx\quad\quad\quad( f\in L^1(G, \omega), h\in L^\infty(G, 1/\omega)).
$$
Let $C_b(G, 1/\omega)$ denote the
subspace of $L^\infty(G, 1/\omega)$ consisting of all bounded continuous
functions, and let $C_0(G, 1/\omega)$ be the
subspace of $C_b(G, 1/\omega)$ consisting of all functions that vanish at
infinity. 
Let also $M(G, \omega)$ be the Banach algebra of all complex regular Borel measures
$\mu$ on $G$ for which $\omega\mu\in M(G)$. It is well-known that
$M(G, \omega)$ is the dual
space of $C_0(G, 1/\omega)$ \cite{dl, r0, sto}, see \cite{r111, rv1} for study of weighted semigroup measure algebras; see also \cite{mr1, mr2, mr}.

 For a Banach algebra $A$, let us recall that a bounded linear operator  $D$ from $A$ into $A^*$ is called a \emph{derivation} if $D(ab)=D(a)\cdot b + a\cdot D(b)$ for all $a, b\in A$. Also, $D$ is said to be an \emph{inner derivation} if there exists $z\in A^*$ such that for every $a\in A$
$$
D(a)=\hbox{ad}_z(a):=z\cdot a-a\cdot z.
$$
The space of all continuous (inner) derivations from $A$ into $A^*$ is denoted by (${\cal I}_{nn}(A, A^*)$) ${\cal Z}( A, A^*)$, respectively. A Banach algebra $A$ is called \emph{weakly amenable} if
$$
{\cal Z}(A, A^*)={\cal I}_{nn}(A, A^*).
$$
Johnson \cite{j} proved that if $G$ is a  locally compact group $G$, then the Banach algebra $L^1(G)$ is weakly amenable; for a simple proof of this result see \cite{dg}. One can arise naturally the question of whether $L^1(G, \omega)$ is weakly amenable. Several authors studied this problem. For example, Bade, Curtis and Dales \cite{bcd} characterized weak amenability of the Banach algebra $L^1({\Bbb Z}, \omega_\alpha)$, where $\omega_\alpha(n)=(1+|n|)^\alpha$. Gronbaek \cite{gro} gave a necessary and sufficient condition for weak amenability of $L^1(G, \omega)$, when $G$ is a discrete Abelian group. He proved that $\ell^1(G, \omega)$ is weakly amenable if and only if there is no non-zero group homomorphism $\phi: G\rightarrow{\Bbb C}$ such that $\phi\in C_b(G, 1/\omega^*)$. The question comes to mind immediately: Dose this result hold for any Abelian group? Recently, Zhang \cite{z} gave
an affirmative answer to this conjecture. It is natural to ask whether this result remain valid for non-Abelian locally compact groups? Borwick \cite{b} studied this question for non-Abelian discrete groups and gave the conditions that characterizes weak amenability of $\ell^1(G, \omega)$; see also \cite{s1} for weak amenability of $\ell^1({\Bbb F}_2, \omega)$ and $\ell^1(\mathbf{(ax+b)}, \omega)$.
In this paper, we answer to this conjecture for non-Abelian locally compact groups.

This paper is organized as follow.  In Section 2, we introduce the notions quasi-additive and inner quasi-additive function. We show that
$${\cal Z}( L^1(G, \omega), L^1(G, \omega)^*)$$
and the space bounded quasi-additive functions are isometrically isomorphic as Banach spaces. This statement holds for $${\cal I}_{nn}(L^1(G, \omega), L^1(G, \omega)^*)$$ and the space of inner quasi-additive functions. Using theses result, we prove that $L^1(G, \omega)$ is weakly amenable if and only if every non-inner quasi-additive function in $L^\infty(G, 1/\omega)$ is unbounded. This cover some known results concerning weak amenability of commutative Beurling algebras. In fact, we give an answer to the question raised in \cite{z}. In Section 3, we consider two Beurling algebras $L^1(G, \omega_i)$, for $i=1, 2$, and investigate the relation between weak amenability of them.
In section 4, we take up the connection between weak amenability of $L^1(G_1, \omega_1)$ and $L^1(G_2, \omega_2)$.

\section{\normalsize\bf Weak amenability of weighted group algebras}

Let $G$ be a locally compact group. A Borel measurable function $p: G\times G\rightarrow{\Bbb C}$ is called \emph{quasi-additive} if for almost every where $x, y, z\in G$
$$p(xy, z)=p(x, yz)+p(y, zx).$$
 A quasi-additive function $p$ is called \emph{inner} if there exists $h\in L^\infty(G, 1/\omega)$ such that $$p(x, y)= h(xy)-h(yx)$$ for almost every where $x, y\in G$.
Let $Q(G)$ be the set of all quasi-additive functions, $D(G, \omega)$ be the set of all $p\in Q(G)$ such that
$$C(p, \omega):=\sup_{x, y\in G}\frac{|p(x, y)|}{\omega^\otimes (x, y)}<\infty$$
and $I(G, \omega)$ be the set of inner quasi-additive functions.
Clearly,
$$
I(G, \omega)\subseteq D(G, \omega)\subseteq Q(G)\cap L^\infty(G\times G, 1/ \omega^\otimes ).
$$
A quasi-additive function $p$ is called \emph{non-inner in} $L^\infty(G, 1/\omega)$ if $p\in Q(G)\setminus I(G, \omega)$.
Let $\check{G}_\omega$ be the set of all group homomorphisms $\phi: G\rightarrow{\Bbb C}$ such that $\phi/\omega$ is bounded. We write $\check{G}$ for $\check{G}_1$.
Note that if  $p\in D(G, \omega)$ , $q\in\check{G}$ and $h\in L^\infty(G)$. Then the functions $ p_1, p_2, p_3: G\times G\rightarrow{\Bbb C}$ defined by $p_1(x, y)=q(x)$,
$$
p_2(x, y)=h(xy)-h(yx)\quad\hbox{and}\quad p_3(x, y)=p(x^{-1}, y^{-1})
$$
are elements in $D(G, \omega)$. Also, if we define the function $q_1: G\rightarrow{\Bbb C}$ by $q_1(x)=p(x, x^{-1})$, then $q_1\in\check{G}$. We now give some properties of quasi-additive functions.

\begin{lemma} \label{zmsj} Let $G$ be a locally compact group and let $\omega$ and $\omega_0$ be weight functions on $G$. Then the following statements hold.

\emph{(i)} $D(G, \omega)$ is a closed subspace of $L^\infty(G\times G, 1/\omega^\otimes )$. In particular, $D(G, \omega)$  is a Banach space.

\emph{(ii)} $I(G, \omega)$ is a subspace of $D(G, \omega)$. Furthermore, if $C=\{h\in L^\infty(G, 1/\omega): h(xy)=h(yx)\;\hbox{for al}l\; x, y\in G\}$, then $I(G, \omega)$ and $L^\infty(G, 1/\omega)/C$ are isomorphic.

\emph{(iii)} If $\omega_0\leq m\omega$ for some $m>0$, then  $D(G, \omega_0)$ and $I(G, \omega_0)$ are subspaces of $D(G, \omega)$ and $I(G, \omega)$, respectively.

\emph{(iv)} $D(G, \omega)$ and $I(G, \omega)$ are subsets of $D(G, \omega^*)$ and $I(G, \omega^*)$, respectively.

\emph{(v)} $\check{G}$ can be embedded into $D(G, \omega)$.
\end{lemma}
{\it Proof.} The statement (i) is proved by using standard arguments. For (ii), we define the function $\Lambda: L^\infty(G, 1/\omega)\rightarrow I(G, \omega)$ by $$\Lambda(h)(x, y)=h(xy)-h(yx).$$ Then $\Lambda$ is an epimorphism with $\hbox{ker}(\Lambda)=C$. Hence
$
L^\infty(G, 1/\omega)/C
$
and $I(G, \omega)$ are isomorphic. Hence (ii) holds. Note that if $\omega_0\leq m\omega$ for some $m>0$, then $C(p, \omega)\leq m^2 C(p, \omega_0)$ and $L^\infty(G, 1/\omega_0)$ is contained in $L^\infty(G, 1/\omega)$. So (iii) holds. The statement (iv) follows from (iii) and the fact that $\omega\leq\omega^*$. Finally, the function $q\mapsto q\circ\pi_1$ from $\check{G}$ onto $D(G, \omega)$ is isomorphism, where $\pi_1$ is the canonical projection. Hence (v) holds.$\hfill\square$\\

Let $G_i$ be a locally compact group  and $\omega_i$ be a weight function on $G_i$ for $i=1, 2$.  We define
$$
\omega_1\otimes\omega_2(x_1, x_2)=\omega_1(x_1)\omega_2(x_2)
$$ for all $x_i\in G_i$. It is easy to prove that $\omega_1\otimes\omega_2$ is a weight function on $G_1\times G_2$.

\begin{lemma} Let $G_i$ be a locally compact group and let $\omega_i$ be a weight functions on $G_i$, for $i=1,2$. Then $D(G_1\times G_2, \omega_1\otimes\omega_2)$ can be embedded into $D(G_1, \omega_1)\times D(G_2, \omega_2)$.
\end{lemma}
{\it Proof.} We only note that the function $p\mapsto (p_1, p_2)$ from $D(G_1\times G_2, \omega_1\otimes\omega_2)$ into  $D(G_1, \omega_1)\times D(G_2, \omega_2)$ is injective, where
$$
p_1(x_1, y_1)=p((x_1, e_2), (y_1, e_2))\quad\hbox{and}\quad p_2(x_2, y_2)=p((e_1, x_2), (e_1, y_2))
$$
for all $x_1, y_1\in G_1$ and $x_2, y_2\in G_2$.$\hfill\square$\\

For Banach algebras $A$ and $B$, let $A\hat{\otimes} B$ be the projective tensor product of $A$ and $B$. Let also ${\cal B}(A, B^*)$ be the space of bounded linear operators from $A$ into $B^*$. Then the function
$$
\Gamma: {\cal B}(A, B^*)\rightarrow (A\hat{\otimes}B)^*
$$
defined by
$$
\langle \Gamma(T), a\otimes b\rangle=\langle T(a), b\rangle
$$
is an isometric isomorphism as Banach spaces; see for example Proposition 13 VI in \cite{bd}. We now give the next theorem which is actually the key to prove our results.

\begin{theorem} \label{zms} Let $G$ be a locally compact group  and $\omega$ be a weight function on $G$. Then the following statements hold.

\emph{(i)} The function $\Gamma: {\cal Z}( L^1(G, \omega), L^\infty(G, 1/\omega))\rightarrow D(G, \omega)$ is an isometric isomorphism, as Banach spaces.

\emph{(ii)} $\Gamma({\cal I}_{nn}(L^1(G, \omega), L^\infty(G, 1/\omega))= I(G, \omega)$.

\emph{(iii)} If $D\in {\cal Z}( L^1(G, \omega), L^\infty(G, 1/\omega))$, then there exists a unique $p\in D(G, \omega)$ such that $D(f)(x)=\int_G p(x, y) f(y) \; dy$ for almost every where $x\in G$.

%\emph{(iv)} If $G$ is amenable and $\omega^*$ is bounded, then $D(G, \omega)=I(G, \omega)$.

%\emph{(v)} If $G$ is Abelian and $\omega^*$ is bounded, then $D(G, \omega)=I(G, \omega)=\{0\}.$

\emph{(iv)}  If  $\omega$  is multiplicative, then $D(G, \omega)=I(G, \omega)$. In particular, $D(G, 1)=I(G, 1).$
\end{theorem}
{\it Proof.} (i)  Let $D\in {\cal Z}( L^1(G, \omega), L^\infty(G, 1/\omega))$. Since $D$ is a bounded linear operator, setting  $A=B=L^1(G, \omega)$ in the definition of $\Gamma$, we have $$p:=\Gamma(D)\in( L^1(G, \omega)\hat{\otimes} L^1(G, \omega))^*=L^\infty(G\times G, 1/\omega^\otimes )$$ and
\begin{eqnarray}\label{main}
\langle D(f), g\rangle=\langle p, f\otimes g\rangle=\int_G\int_G p(x, y) f(x)g(y) dxdy
\end{eqnarray}
for all $f, g\in L^1(G, \omega)$. On the other hand, for every $f,g, k\in L^1(G, \omega)$ we have
\begin{eqnarray*}
\langle D(f\ast g), k\rangle&=&\langle D(f)\cdot g, k\rangle+\langle  f\cdot D(g), k\rangle\\
&=&\langle D(f), g\ast k+k\ast f\rangle.
\end{eqnarray*}
It follows that
$$
\int_G\int_G p(x, y) (f\ast g)(x) k(y)\;dxdy=\int_G\int_G p(x, y)  f(x) (g\ast k+ k\ast f)(y)\; dxdy
$$
Thus for every $f,g,h\in L^1(G, \omega)$, we obtain
\begin{eqnarray*}
\int_G\int_G\int_G p(xy, z) f(x)g(y)k(z)\;dxdydz&=&\int_G\int_G\int_G p(x,yz)f(x)g(y)k(z)\;dxdydz\\
&+& \int_G\int_G\int_G p( y,zx)f(x)g(y)k(z)\;dxdydz.
\end{eqnarray*}
This implies that $p\in D(G, \omega)$. Therefore, $\Gamma$ maps $ {\cal Z}( L^1(G, \omega), L^\infty(G, 1/\omega))$ into $D(G, \omega)$, as Banach spaces.

Now, let $p\in D(G, \omega)$. Define the linear operator $D: L^1(G, \omega)\rightarrow L^\infty(G, 1/\omega)$ by
$$
\langle D(f), g\rangle:=\langle p, f\otimes g\rangle
$$
for all $f, g\in L^1(G, \omega)$. Since $p\in D(G, \omega)$, a similar argument as given above follows that  $D\in {\cal Z}( L^1(G, \omega), L^\infty(G, 1/\omega))$. So $\Gamma(D)=p$. Thus $\Gamma$ is an isometric isomorphism from ${\cal Z}( L^1(G, \omega), L^\infty(G, 1/\omega))$ onto $D(G, \omega)$.

(ii) First, note that $L^\infty(G, 1/\omega)$ is a Banach $L^1(G, \omega)-$bimodule with the following actions.
$$
f\cdot h(x)=\int_Gf(y) h(xy)\; dy\quad\hbox{and}\quad h\cdot f(x)=\int_G f(y) h(yx)\;dy
$$
for all $f\in L^1(G, \omega)$, $h\in L^\infty(G, 1/\omega)$ and $x\in G$.

Assume note that $D\in {\cal I}_{nn}(L^1(G, \omega), L^\infty(G, 1/\omega)$. Then $D=\hbox{ad}_h$ for some $h\in L^\infty(G, 1/\omega)$.
So
\begin{eqnarray}\label{x}
\langle D(f), g\rangle&=&\langle h\cdot f-f\cdot h, g\rangle\nonumber\\
&=&\int_G g(y) (h\cdot f(x)-f\cdot h(x))\; dx\\
&=&\int_G\int_G g(y) f(x)(h(xy)-h(yx))\; dx dy\nonumber
\end{eqnarray}
for all $f, g\in L^1(G, \omega)$.
Consequently, $\langle D(f), g\rangle=\langle p, f\otimes g\rangle$, where $$p(x, y)=h(xy)-h(yx)$$ for all $x, y\in G$. Obviously, $p\in D(G, \omega)$. This implies that  $$\Gamma({\cal I}_{nn}(L^1(G, \omega), L^\infty(G, 1/\omega))= I(G, \omega).$$

(iii) Let $D\in {\cal Z}( L^1(G, \omega), L^\infty(G, 1/\omega))$. If $f\in L^1(G, \omega)$, then $D(f)\in L^1(G, \omega)^*$. So for every $g\in L^1(G, \omega)$, we have
$$
\langle D(f), g\rangle=\int_G D(f)(y)g(y)\; dy.
$$
From this and (\ref{main}) we conclude that (iii) holds.

%(iv) Let $G$ be amenable and $\omega^*$ be bounded. Then $L^1(G, \omega)$ is amenable and so it is weakly amenable. Hence
%$$
%{\cal Z}( L^1(G, \omega), L^\infty(G, 1/\omega))={\cal I}_{nn}(L^1(G, \omega), L^\infty(G, 1/\omega)).
%$$
%Thus $D(G, \omega)=I(G, \omega)$.

%(v) If $G$ is Abelian, then ${\cal I}_{nn}(L^1(G, \omega), L^\infty(G, 1/\omega))=\{0\}$. This together with (iv) follows (v).

(iv) Let $\omega$ be multiplicative. Then $L^1(G, \omega)$ and $L^1(G)$ are isomorphism, as Banach algebras. Since $L^1(G)$ is weakly amenable, the statement (v) holds.
$\hfill\square$\\

We now state the main result of this section which answers to an open problem given in \cite{z}.

%\begin{theorem}\label{xx2} Let $G$ be a locally compact group and $\omega$ be a weight function on $G$. Then the following assertions are equivalent.

%\emph{(a)} $L^1(G, \omega)$ is weakly amenable.

%\emph{(b)} For every $D\in {\cal Z}(  L^1(G, \omega), L^\infty(G, 1/\omega))$ there exists $h\in L^\infty(G, 1/\omega)$ such that for every $f, g\in L^1(G, \omega)$
% $$\langle D(f), g\rangle=\int_G\int_G f(x)g(y)(h(xy)-h(yx))\;dxdy.$$
%\quad~\emph{(c)} $D(G, \omega)= I(G, \omega)$.

%\emph{(d)} $\{p\in Q(G): C(p,\omega)<\infty\}=I(G, \omega)$.
%\end{theorem}
%{\it Proof.}
\begin{theorem}\label{xx2} Let $G$ be a locally compact infinite group and $\omega$ be a weight function on $G$. Then the following assertions are equivalent.

\emph{(a)} $L^1(G, \omega)$ is weakly amenable.

\emph{(b)} For every bounded derivation $D: L^1(G, \omega)\rightarrow L^\infty(G, 1/\omega)$ there exists $h\in L^\infty(G, 1/\omega)$ such that $\langle D(f), g\rangle=\int_G\int_G f(x)g(y)(h(xy)-h(yx))\;dxdy$ for all $f, g\in L^1(G, \omega)$.

\emph{(c)} For every quasi-additive function $p\in D(G, \omega)$ there exists $h\in L^\infty(G, 1/\omega)$ such that $p(x, y)= h(xy)-h(yx)$ for all $x, y\in G$.

\emph{(d)} Every quasi-additive function $p$ with $C(p, \omega)<\infty$ is inner.

\emph{(e)} Every non-inner quasi-additive function in $L^\infty(G, 1/\omega)$ is unbounded.

\emph{(f)} $D(G, \omega)=I(G, \omega)$.
\end{theorem}
{\it Proof.} Let $D\in {\cal Z}(  L^1(G, \omega), L^\infty(G, 1/\omega))$. If  $L^1(G, \omega)$ is weakly amenable, then $D\in {\cal I}_{nn}(L^1(G, \omega), L^\infty(G, 1/\omega))$. By Theorem \ref{zms} (i), $\Gamma(D)\in I(G, \omega)$. Thus there exists $h\in L^\infty(G, 1/\omega)$ such that
$$
\Gamma(D)(x, y)=h(xy)-h(yx)
$$
for all $x, y\in G$. It follows from (\ref{main}) that
$$
\langle D(f), g\rangle=\int_G\int_G f(x)g(y) (h(xy)-h(yx))\; dxdy
$$
for all $f, g\in L^1(G, \omega)$. So (a)$\Rightarrow$(b).

Let $p\in D(G, \omega)$. Then $\Gamma(D)=p$ for some $D\in {\cal Z}(  L^1(G, \omega), L^\infty(G, 1/\omega))$. If (b) holds, then by (\ref{main}), for every $f, g\in L^1(G, \omega)$, we have
\begin{eqnarray*}
\int_G\int_G p(x, y) f(x) g(y)\; dx dy&=&\langle\Gamma(D), f\otimes g\rangle=\langle D(f), g\rangle\\
&=&\int_G\int_G f(x)g(y)(h(xy)-h(yx))\;dxdy.
\end{eqnarray*}
This shows that $p(x, y)= h(xy)-h(yx)$ for almost every where $x, y\in G$. So (b)$\Rightarrow$(c).
For the implication (c)$\Rightarrow$(d), note that
$$
D(G,\omega)=\{p\in Q(G): C(p,\omega)<\infty\}.
$$
The implications (d)$\Rightarrow$(e)$\Rightarrow$(f) are clear. Finally, let $D\in {\cal Z}(  L^1(G, \omega), L^\infty(G, 1/\omega))$. Then $\Gamma(D)=p$ for some $p\in D(G, \omega)$. If $p\in I(G, \omega)$, then $D\in {\cal I}_{nn}( L^1(G, \omega), L^\infty(G, 1/\omega))$. That is, (d) implies (a).
$\hfill\square$\\

As a consequence of Theorem \ref{xx2} we give the following result.

\begin{corollary} Let $G$ be a locally compact infinite group. Then $L^1(G, \omega)$ is not weakly amenable if and only if there exists a bounded non-inner quasi-additive function in $L^\infty(G, 1/\omega)$.
\end{corollary}

The next result can be considered as an improvement of Theorem 3.1 in \cite{z}.

\begin{corollary}\label{ab} Let $G$ be a locally compact Abelian group. Then the following assertions are equivalent.

\emph{(a)} $L^1(G, \omega)$ is weakly amenable.

\emph{(b)} The zero map is the only quasi-additive function in $D(G, \omega)$.

\emph{(c)} The zero map is the only group homomorphism in $L^\infty(G, 1/\omega^*)$.

\emph{(d)} The zero map is the only group homomorphism in $C_b(G, 1/\omega^*)$.
\end{corollary}
{\it Proof.} From Theorem  1.2 in \cite{z} and Theorem \ref{xx2}, the implications (d)$\Rightarrow$(a)$\Rightarrow$(b) hold. The implication (c)$\Rightarrow$(d) is clear. To show that (b)$\Rightarrow$(c), let $\phi$ be a group homomorphism in $L^\infty(G, 1/\omega^*)$. Define $p: G\times G\rightarrow{\Bbb C}$ by $p(x, y)=\phi(x)$. Then $p$ is a quasi-additive function in $L^\infty(G\times G, 1/\omega^\otimes )$. So $\phi=0$.$\hfill\square$\\

In the sequel, we give another consequent of Theorem \ref{xx2}.

\begin{corollary}  Let $G$ be a locally compact infinite group. Then the following statements hold.

\emph{(i)} Every non-inner quasi-additive function in $L^\infty(G)$ is unbounded.

\emph{(ii)} If $G$ is Abelian, then every non-zero group homomorphism from $G$ into ${\Bbb C}$ is unbounded.
\end{corollary}
{\it Proof.} It is well-known that if $G$ is a locally compact group, then $L^1(G)$ is weakly amenable; see for example \cite{dl}. This together with Theorem \ref{xx2} establish (i). The statement (ii) follows from Corollary \ref{ab}.$\hfill\square$

\begin{example}{\rm Define the function $\omega: \Bbb{R}\rightarrow [1, \infty) $ by
\begin{eqnarray*}
\omega(x)=\left\{
\begin{array}{rl}
\hbox{exp}(x) & x\geq 0\\
1    & x<0
\end{array}\right.
\end{eqnarray*}
Then $\omega$ is a weight function on $\Bbb{R}$. For $r\in\Bbb{R}$, we define $p_r(x)=rx$ for all $x\in\Bbb{R}$. Then
$$
\sup{\frac{|p_r(x)|}{\omega^*(x)}: x\in\Bbb{R}\}=\sup\{|rx|}{\hbox{exp}(x)}: x\in\Bbb{R}^+\}\leq|r|.
$$
Hence there exists a non-zero group homomorphism in $L^\infty(\Bbb{R}, 1/\omega^*)$. Thus $L^1(\Bbb{R}, \omega)$ is not weakly amenable. It is proved that $$D(L^1(G, \omega))=[-1,0]\times\Bbb{R}.$$
}
\end{example}

\section{\normalsize\bf Beurling algebras with different weights}\label{sec4}

We commence this section with the following result which is interesting and useful.

\begin{theorem}\label{con} Let $\omega_1$ and $\omega_2$ be weight functions on a locally compact group $G$. Then the following assertions are equivalent.

\emph{(a)} $L^1(G, \omega_2)$ is a subspace of $L^1(G, \omega_1)$.

\emph{(b)} $M(G, \omega_2)$ is a subspace of $M(G, \omega_1)$.

\emph{(c)} There exists $m>0$ such that $\|f\|_{\omega_1}\leq m\|f\|_{\omega_2}$ for all $f\in L^1(G, \omega_2)$.

\emph{(d)} There exists $n>0$ such that $\|\mu\|_{\omega_2}\leq n\|\mu\|_{\omega_1}$ for all $\mu\in M(G, \omega_2)$.

\emph{(e)} There exists $s>0$ such that $\omega_1\leq s\omega_2$.
\end{theorem}
{\it Proof.} The implications (e)$\Rightarrow$(c)$\Rightarrow$(a) are clear. If $L^1(G, \omega_2)$ is a subspace of $L^1(G, \omega_1)$, we can regard $L^1(G, \omega_2)^{**}$ as a subspace of $L^1(G, \omega_1)^{**}$. Hence $M(G, \omega_1)\oplus C_0(G, 1/\omega_1)^\perp$ is contained in $M(G, \omega_2)\oplus C_0(G, 1/\omega_2)^\perp$. This shows that $M(G, \omega_1)=
\pi_1(L^1(G, \omega_1)^{**})$ is a subset of $M(G,\omega_2)=\pi_1(L^1(G, \omega_2)^{**})$. That is, (a)$\Rightarrow$(b). To show that (b)$\Rightarrow$(d), we define $$\|\mu\|=\|\mu\|_{\omega_1}+\|\mu\|_{\omega_2},$$ then $B:=(M(G, \omega_2), \|.\|)$ is a Banach algebra. So the identity map $F: B\rightarrow A$ is injective and continuous, where $A:= (M(G, \omega_1), \|.\|_{\omega_1})$. So $F: B\rightarrow F(B)$ is bijective. In view of the inverse mapping theorem, $F^{-1}$ is continuous. Thus there exists $r>1$ such that
$$
\|F^{-1}(\mu)\|\leq r\|\mu\|_{\omega_1}
$$
for all $\mu\in M(G, \omega_1)$. Consequently,
$$
\|\mu\|_{\omega_1}+\|\mu\|_{\omega_2}\leq r\|\mu\|_{\omega_1}.
$$
Therefore, $\|\mu\|_{\omega_2}\leq (r-1)\|\mu\|_{\omega_1}$.
%   For every $f\in L^1(G, \omega_2)$, we define
%$$
%\||f \||=\|f\|_{\omega_1}+\|f\|_{\omega_2}.
%$$
%Then the identity map
%$I:(L^1(G,\omega_2),\||.\||)\rightarrow L^1(G, \omega_2)$ is continuous and bijective. In view of the inverse mapping theorem, $I^{-1}$ is continuous.
%Hence there exists $M>1$ such that $
%\||f\||\leq M\|f\|_{\omega_2}
%$ for all $f\in L^1(G, \omega_2)$. Therefore, $\|f\|_{\omega_1}\leq(M-1)\|f\|_{\omega_2}$. That is (a) implies (c). Similarly, (b) implies (d).
Finally, let (d) hold. Then for every $x\in G$, we have $\|\delta_x\|_{\omega_2}\leq n\|\delta_x\|_{\omega_1}$. Thus $\omega_1\leq n\omega_2$ for some $n>0$. That is, (d)$\Rightarrow$(e).$\hfill\square$\\

Let us recall that two weight function $\omega_1$ and $\omega_2$ are \emph{equivalent} if $m\omega_1\leq\omega_2\leq n\omega_1$ for some $m, n>0$.

\begin{corollary} Let $\omega_1$ and $\omega_2$ be  weight functions on a locally compact group $G$. Then the following assertions are equivalent.

\emph{(a)} $L^1(G, \omega_1)=L^1(G, \omega_2)$.

\emph{(b)} $\|.\|_{\omega_1}$ and $\|.\|_{\omega_2}$ are equivalent.

\emph{(c)} $\omega_1$ and $\omega_2$ are equivalent.
\end{corollary}

In the following, we consider two weight functions $\omega_1$ and $\omega_2$ on a locally compact group $G$ and study the relation between weak amenability of $L^1(G, \omega_1)$ and $L^1(G, \omega_2)$.

\begin{proposition}\label{w**} Let $\omega_1$ and $\omega_2$ be weight functions on a locally compact group $G$ such that $\omega_1\leq m\omega_2$ for some $m>0$. Then the following statements hold.

\emph{(i)} If $L^1(G, \omega_2)$ is weakly amenable and $I(G, \omega_1)=I(G, \omega_2)$, then $L^1(G, \omega_1)$ is weakly amenable.

\emph{(ii)} If $\omega_2$ or $\omega_2^*$ is bounded, then $L^1(G, \omega_1)$ is weakly amenable.
\end{proposition}
{\it Proof.} (i) By Lemma \ref{zmsj} and Theorem \ref{xx2} we have
$$
D(G, \omega_1)\subseteq D(G, \omega_2)=I(G, \omega_2)=I(G, \omega_1)\subseteq D(G, \omega_1).
$$
Hence $D(G, \omega_1)= I(G, \omega_1)$ and so $L^1(G, \omega_1)$ is weakly amenable.

(ii) Let $\omega_2$ be bounded. Then $\omega_1$ is bounded and so $L^1(G, \omega_1)$ is weakly amenable. Now if $\omega_2^*$ is bounded, then $\omega_2$ is bounded, because $\omega_2\leq\omega_2^*$. Hence $L^1(G, \omega_1)$ is weakly amenable.$\hfill\square$\\

As an immediate corollary of  Proposition \ref{w**} we have the following result.

\begin{corollary} Let $\omega_1$ and $\omega_2$ be weight functions on a locally compact group $G$ such that $\omega_1$ and $\omega_2$ are equivalent. Then  $L^1(G, \omega_1)$ is weakly amenable if and only if $L^1(G, \omega_2)$ is weakly amenable.
\end{corollary}

\begin{example}{\rm Let $\omega$ be a weight function on a locally compact group $G$ and $a\in G$. The function $\omega_a$ defined by $\omega_a(x)=\omega(axa^{-1})$ is a weight function on $G$. For every $x\in G$, we have
\begin{eqnarray*}
\omega(x)=\omega(a^{-1}(axa^{-1})a)\leq\omega(a^{-1})\omega_a(x)\omega(a)
\end{eqnarray*}
and
$$
\omega_a(x)=\omega(axa^{-1})\leq\omega(a)\omega(a^{-1})\omega(x).
$$
Hence
\begin{eqnarray}\label{12345}
(1/m)\omega(x)\leq\omega_a(x)\leq m\omega(x),
\end{eqnarray}
where $m=\omega(a)\omega(a^{-1})$. That is, $\omega$ and $\omega_a$ are equivalent. Now, define the function $W$ on $G$ by $$W(x)=\inf_{a\in G} \omega_a(x).$$ Then $W$ is a weight function on $G$. By (\ref{12345}), $\omega$ and $W$ are equivalent. Hence $L^1(G, \omega)$ is weakly amenable if and only if $L^1(G, \omega_a)$ is weakly amenable; or equivalently, $L^1(G, W)$ is weakly amenable.}
\end{example}

As another consequence of Theorem \ref{w*} we have the following result. The part (i) due to Pourabbas \cite{p}.

\begin{corollary}\label{w*} Let $G$ be a locally compact group. Then the following statements hold.

\emph{(i)} If $\omega^*$ is bounded, then  $L^1(G, \omega)$ is weakly amenable and $D(G, \omega^*)=I(G, \omega^*)$.

\emph{(ii)} If $L^1(G, \omega^*)$ is weakly amenable and $I(G, \omega)=I(G, \omega^*)$, then $L^1(G, \omega)$ is weakly amenable and $D(G, \omega)=D(G, \omega^*)$. In the case where $G$ is Abelian, $D(G, \omega)=\{0\}$.

\emph{(iii)} If $D(G, \omega^*)=I(G, \omega^*)=I(G, \omega)$, then $D(G, \omega)=I(G, \omega)$.
\end{corollary}
{\it Proof.}  %Let $p\in Q(G)\setminus I(G, \omega)$. Set
%$$
%{\cal C}=\{\frac{|p(x,y)|}{\omega^\otimes(x, y)\omega^\otimes(a^{-1}, b^{-1})}: a, b, x, y\in G\}
%$$
%and
%$$
%{\cal D}=\{\frac{|p(a^{-1},b^{-1})|}{\omega^\otimes( x^{-1}, y^{-1})\omega^\otimes(a, b)}: a, b, x, y\in G\}.
%$$
% If $L^1(G, \omega^*)$ is weakly amenable, then
%$$
%\infty=\sup({\cal CD})=\sup({\cal C})\sup({\cal D})= \sup({\cal C})^2.
%$$
%So $\sup(C)=\infty$ and thus $\sup(A)=\infty$. This shows that $L^1(G,\omega)$ is weakly amenable.
Since $\omega\leq\omega^*$, the statement (i) follows from Proposition \ref{w**} (ii). From Proposition \ref{w**} yields that $L^1(G, \omega)$ is weakly amenable and by Theorem \ref{xx2},
$$
D(G, \omega)=I(G, \omega)=I(G, \omega^*)=D(G, \omega^*).
$$
When $G$ is Abelian, we have
$$
{\cal Z}( L^1(G, \omega), L^\infty(G, 1/\omega))= {\cal I}_{nn}(L^1(G, \omega), L^\infty(G, 1/\omega))=\{0\}.
$$
In view of Theorem \ref{zms}, we infer that $D(G, \omega)=I(G, \omega)=\{0\}$. So (ii) holds. For (iii), by Lemma \ref{zmsj} we have
$$
D(G, \omega)\subseteq D(G, \omega^*)=I(G, \omega^*)=I(G, \omega)\subseteq D(G, \omega).
$$
Hence (iii) holds.$\hfill\square$

\begin{example}{\rm It is well-known from \cite{bcd} that $\ell^1({\Bbb Z}, \omega_\alpha)$ is weakly amenable if and only if $0\leq\alpha<1/2$. Hence the Banach algebra $\ell^1({\Bbb Z}, \omega_{1/3})$ is weakly amenable, however, $\ell^1({\Bbb Z}, \omega_{1/3}^*)$  is not weakly amenable. Also, $\ell^1({\Bbb Z}, \omega_{1/6}^*)$  is weakly amenable, but, $\omega_\alpha^*$ is unbounded.}
\end{example}

Let $\omega$ be a weight function on $G$. Then $\omega^\prime$ defined by $\omega^\prime(x)=\omega(x^{-1})$ for all $x\in G$ is also a weight function on $G$.

\begin{theorem}\label{d} Let $G$ be a locally compact group and $\omega$ be a weight function on $G$. Then $L^1(G, \omega^\prime)$ is weakly amenable if and only if $L^1(G, \omega)$ is weakly amenable.
\end{theorem}
{\it Proof.} Let $p$ be a complex-valued function $G\times G$. Define $p^\prime(x, y)=p(x^{-1}, y^{-1})$ for all $x, y\in G$. It routine to check that
$$C(p, \omega)= C(p^\prime, \omega^\prime).$$
Now, for a complex-valued function $h$ on $G$, define $h^\prime(x)=-h(x^{-1})$ for all $x\in G$. Obviously,
$p(x, y)=h(xy)-h(yx)$ if and only if $h^\prime(x, y)=h^\prime(xy)-h^\prime(yx)$ for all $x, y\in G$. Therefore, $p\in I(G, \omega)$ if and only if $p^\prime\in I(G, \omega^\prime)$. These facts together with Theorem \ref{xx2} prove the result.$\hfill\square$\\

It is easy to see that if $G$ is a locally compact Abelian group, then $\frak{D}:=\{(x, x^{-1}): x\in G\}$ is a locally compact Abelian group.

\begin{theorem}\label{d} Let $G$ be a locally compact Abelian group and $\omega$ be a weight function on $G$. Then $L^1(G, \omega^*)$ is weakly amenable if and only if $L^1(\frak{D}, \omega^\otimes )$ is weakly amenable.

\emph{(ii)} $L^1(G, \omega^\prime)$ is weakly amenable if and only if $L^1(G, \omega)$ is weakly amenable, where $\omega^\prime(x)=\omega(x^{-1})$ for all $x\in G$.
\end{theorem}
{\it Proof.} Let $q: G\rightarrow{\Bbb C}$ be a group homomorphism. Define the function $Q: \frak{D}\rightarrow{\Bbb C}$ by $Q(x, x^{-1})=q(x)$. Then $Q$ is a group homomorphism. Since for every $x\in G$
$$
\frac{|q(x)|}{(\omega^*)^*(x)}= \frac{| Q(x,x^{-1})|}{(\omega^\otimes)^* (x,x^{-1})},
$$
it follows that $q\in\check{G}_{(\omega^*)^*}$ if and only if $Q\in\check{\frak{D}}_{(\omega^\otimes)^*}$. Also, $q=0$ if and only if $Q=0$. Now, apply Corollary \ref{ab}.
%, $Q=0$ and so $q=0$. Thus $L^1(G, \omega^*)$ is weakly amenable.
%To converse, suppose that $L^1(\frak{D}, \omega^\otimes )$  is not weakly amenable. Since $\frak{D}$ is Abelian, there exists a non-zero element $Q\in\check{\frak{D}}$ such that
%$$
%\sup{|Q( x,x^{-1})|/(\omega\times\omega)^*}<\infty.
%$$
%Define $q(x):=Q( x,x^{-1})$ for all $x\in G$. Then $q\in\check{G}$ and
%$\sup{|q(x)|/\omega^{**}(x)}<\infty.$ It follows that $L^1(G, \omega^*)$ is not weakly amenable, a contradiction.
$\hfill\square$

\section{\normalsize\bf Weak amenability of Beurling algebras}\label{sec5}

Let $\phi: G_1\rightarrow G_2$ be a group epimorphism and $\omega$ be a weight function on $G_2$. Then the function $\overleftarrow{\omega}: G_1\rightarrow [1, \infty)$
defined by $\overleftarrow{\omega}(x_1)=\omega(\phi(x_1))$ is a weight function on $G_1$. Define the function
$\frak{S}: Q(G_2)\rightarrow Q(G_1)$ by
$$
\frak{S}(p_2)(x_1, y_1)=p_2(\phi(x_1), \phi(y_1)).
$$
It is clear that $\frak{S}$ is injective and $C(p_2, \omega)=C(\frak{S}(p_2), \overleftarrow{\omega})$ for all $p_2\in Q(G_2)$. So $\frak{S}$ maps $D(G_2, \omega)$ into $D(G_1, \overleftarrow{\omega})$. Since $h_2\circ\phi\in L^\infty(G_1, 1/\overleftarrow{\omega})$ for all $h_2\in L^\infty(G_2, 1/\omega)$. This shows that $\frak{S}(I(G_2, \omega))$ is contained in $I(G_1, \overleftarrow{\omega})$.

\begin{proposition}\label{o} Let $G_1$ and $G_2$ be locally compact infinite groups and let $\omega$ be a weight function on $G_2$. Then the following statements hold.

\emph{(i)} If $L^1(G_2, \omega)$ is weakly amenable, then $L^1(G_1, \overleftarrow{\omega})$ is weakly amenable.

\emph{(ii)} $L^1(G_1, \overleftarrow{\omega})$ is weakly amenable and $\frak{S}(I(G_2, \omega))= I(G_1, \overleftarrow{\omega})$, then $L^1(G_2, \omega)$ is weakly amenable.
\end{proposition}
{\it Proof.} (i) Let $p_1$ be a non-inner quasi-additive function in $L^\infty(G_1, 1/\overleftarrow{\omega})$. Since $\frak{S}$ maps $I(G_2, \omega)$ into $I(G_1, \overleftarrow{\omega})$, there exists a non-inner quasi-additive function $p_2$ in $L^\infty(G_2, 1/\omega)$ such that $\frak{S}(p_2)=p_1$. So if $L^1(G_2, \omega)$ is weakly amenable, then
$$
C(p_1, \overleftarrow{\omega})=C(\frak{S}(p_2), \overleftarrow{\omega})=C(p_2, \omega)=\infty.
$$
Therefore, $L^1(G_1, \overleftarrow{\omega})$ is weakly amenable.

(ii) Let $p_2$ be a non-inner quasi-additive function in $L^\infty(G_2, 1/\omega)$. Since $\frak{S}$ is injective, $\frak{S}(p_2)$ is a non-inner quasi-additive function in $L^\infty(G_1, 1/\overleftarrow{\omega})$. Hence if $L^1(G_1, \overleftarrow{\omega})$ is weakly amenable, then
$$
C(p_2, \omega)=C(\frak{S}(p_2), \overleftarrow{\omega})=\infty.
$$
This shows that $L^1(G_2, \omega)$ is weakly amenable.$\hfill\square$\\

Let $G$ be a locally compact group and $N$ be a normal subgroup of $G$. It is easy to see that the function $$\hat{\omega}(xN)=\inf\{\omega(y): yN=xN\}$$ is a weight function on $G/N$. Also the function $\overline{\omega}(x)=\hat{\omega}(xN)$ is a weight function on $G$ and $\overline{\omega}(x)\leq\omega(x)$ for all $x\in G$.

\begin{corollary} Let $\omega$ be a weight function on a locally compact group $G$ and $N$ be a normal subgroup of $G$ such that $G/N$ is infinite. If $L^1(G/N, \hat{\omega})$ is weakly amenable, then $L^1(G, \overline{\omega})$ is weakly amenable.
\end{corollary}
{\it Proof.} Let $\pi: G\rightarrow G/N$ be the quotient map. Then  for every $x\in G$, we have
$$
\overleftarrow{\omega}(x)= \hat{\omega}(\pi(x))=\hat{\omega}(xN)=\overline{\omega}(x).
$$
Hence $\overleftarrow{\omega}= \overline{\omega}$. Now, invoke Proposition \ref{o}.$\hfill\square$\\

Let $\omega_i$ be a weight function on a locally compact group $G_i$, for $i=1, 2$. Then $W_i(x_1, x_2)=\omega_i(x_i)$ is a weight function on $G_1\times G_2$, for $i=1,2$.

\begin{corollary} Let $\omega_i$ be a weight function on a locally compact infinite group $G_i$, for $i=1, 2$. If $L^1(G_1, \omega_1)$ or $L^1(G_2, \omega_2)$ is weakly amenable, then $L^1(G_1\times G_2, W_i)$ is weakly amenable.
\end{corollary}
{\it Proof.} Let $\pi_i: G_1\times G_2\rightarrow G_i$ be the canonical projection, for $i=1,2$. Then
$$
\overleftarrow{\omega}(x_1, x_2)=\omega_i(\pi_i(x_1, x_2))=\omega_i(x_i)=W_i(x_1, x_2)
$$
for all $x_i\in G_i$. So $\overleftarrow{\omega}=W_i$. By Proposition \ref{o}, the result hold.$\hfill\square$\\

%\begin{corollary} Let $\omega$ be a weight function on a locally compact group $G$ and $N$ be a normal subgroup of $G$. If $L^1(G, \omega)$ is weakly amenable, then $L^1(G/N, \hat{\omega})$ is weakly amenable.
%\end{corollary}
%{\it Proof.}  Let $P\in Q(G/N)$. Then for every $x, y\in G$, we have
%$$
%\frac{|P(xH,yH)|}{(\omega^\prime)^\otimes( xH, yH)}\geq
%\frac{|p(x,y)|}{\omega^\otimes(x, y)},
%$$
%where
%$$
%p( x,y)=P(xH,yH).$$ Note that $p\in Q(G)$ and
%$
%C(p,\omega)=\infty$. Thus
%$
%C(P, \omega^\prime)=\infty.$ Therefore,
%$L^1(G/H, \omega^\prime)$ is weakly amenable. So (i) holds.

%\begin{proposition} Let $G$ be a locally compact group, $H$ be a subgroup of $G$ and $\omega$ be a weight function on $G$. Then the following statements hold.

%\emph{(i)} $D(G,w)$ and $I(G,w)$ can be embedded into $D(H, w|_H)$ and $I(H, \omega|H)$, respectively.

%\emph{(ii)} $C(p|H, \omega|_H)\leq C(p, \omega)$ for all $p\in Q(G)$.

%\emph{(iii)} For every $p\in D(G, \omega)$, $C(p|_{H\times H},\omega|_H)=\infty$, then $L^1(G, \omega)$ is weakly amenable.
%\end{proposition}

%\begin{corollary} Let $\omega$ be a weight function on a locally compact group $G$. Then $L^1(G, \omega)$ is weakly amenable if and only if $L^1(H, \omega|_H)$ is weakly amenable for all subgroup $H$ of $G$.
%\end{corollary}
In the following we investigate the relation between weak amenability of $L^1(G, \omega)$ and $L^1(H, \omega|_H)$ when $H$ is a subgroup of $G$. To this end, we need some results.

\begin{proposition}\label{ex} Let $\omega$ be a weight function on a locally compact Abelian group $G$ and $H$ be a subgroup of $G$. Then every group homomorphism $q: H\rightarrow{\Bbb C}$ has an extension to a group homomorphism $Q: G\rightarrow{\Bbb C}$.
\end{proposition}
{\it Proof.} Let $q: H\rightarrow{\Bbb C}$ be a group homomorphism. Let $A$ be a subset of $G$ such that $\{ Ha: a\in A\}$ is a family of pairwise disjoint subsets of $G$ and $G=\cup_{a\in A} Ha$. If $x\in G$, then there exists a unique element $a_x\in A$ such that $x\in Ha_x$. So there exists a unique element $\tilde{x}=x(a_x)^{-1}$ in $H$ such that $x=\tilde{x} a_x$. Define $Q(x)=q(\tilde{x})$. For every $x, y\in G$, we have $Ha_x Ha_y= Ha_xa_y$. Thus $xy\in  Ha_xa_y\cap Ha_{xy}$. Hence $ Ha_xa_y= Ha_{xy}$. It follows that $a_xa_y= a_{xy}$. Hence
$$\tilde{x} \tilde{y} a_xa_y=
\tilde{x} a_x\tilde{y} a_y=
xy=
\widetilde{xy} a_{xy}=
\widetilde{xy} a_xa_y.
$$
and thus $\widetilde{xy}=\tilde{x} \tilde{y}$. Consequently,
$$Q(xy)=
q(\widetilde{xy})=
q(\tilde{x})+q(\tilde{y})=
Q( x)+ Q(y).$$
This shows that $Q$ is a group homomorphism. For every $h\in H$, we have $he=h$. Thus $a_x=e$. Therefore, $x=\tilde{x}$. That is, $Q|_H=q$.$\hfill\square$\\

Let $G$ be a locally compact Abelian group and $H$ be a subgroup of $G$. We denote by $\widetilde{H}$ be the set of all $\tilde{x}$ pointed out in the proof of Proposition \ref{ex}.

\begin{lemma} Let $G$ be a locally compact Abelian group and $\omega$ be a weight function on $G$. Then the following statements hold.

\emph{(i)} $\widetilde{H}$ is a normal subgroup of $H$.

\emph{(ii)} $\omega^*(\tilde{x})=\omega\otimes\omega(\tilde{x}, (\tilde{x})^{-1})$ for all $x\in G$.

\emph{(iii)} The function $\tilde{\omega}$ defined by $\tilde{\omega}(x)=\omega(\tilde{x})$ is a weight function on $G$. Furthermore, $(\tilde{\omega})^*(x)=\omega^*(\tilde{x})$ for all $x\in G$.
\end{lemma}
{\it Proof.} It is easy to see that $\tilde{e}=e$ and $\widetilde{xy}=\tilde{x} \tilde{y}$ for all $x, y\in G$. since $xx^{-1}=e$, we have $\tilde{x}\widetilde{x^{-1}}=\tilde{e}=e$. Thus
$\widetilde{x^{-1}}=( \tilde{x})^{-1}$. These facts show that $\widetilde{H}$ is a subgroup of $G$. So (i) holds.  It is easy to prove that (ii) holds. Since $\widetilde{xy}=\tilde{x} \tilde{y}$, it follows that
\begin{eqnarray*}
\tilde{\omega}(xy)&=&\omega(\widetilde{xy})=\omega(\tilde{x} \tilde{y})\\
&\leq&\omega(\tilde{x})\omega(\tilde{y})= \tilde{\omega}(x)\tilde{\omega}(y).
\end{eqnarray*}
This shows that $\tilde{\omega}$ is a weight function on $G$. For every $x\in G$, we have
\begin{eqnarray*}
(\tilde{\omega})^*(x)&=&\tilde{\omega}(x) \tilde{\omega}(x^{-1})\\
&=&\omega(\tilde{x})\omega(\widetilde{x^{-1}})\\
&=&\omega(\tilde{x})\omega((\tilde{x})^{-1})\\
&=&\omega^*(\tilde{x}).
\end{eqnarray*}
Therefore, $(\tilde{\omega})^*(x)=\omega^*(\tilde{x})$ for all $x\in G$. That is, (iii) holds. $\hfill\square$

\begin{theorem}\label{til} Let $G$ be a locally compact Abelian group and $H$ be a subgroup of $G$. If $L^1(G, \tilde{\omega})$ is weakly amenable, then $L^1(H, \omega|_H)$ is weakly amenable.
\end{theorem}
{\it Proof.} Let $q: H\rightarrow \Bbb {C}$ be a non-zero group homomorphism. By Proposition \ref{ex}, there exists $Q: G\rightarrow \Bbb {C}$ such that $Q(x)=q(\tilde{x})$ for all $x\in G$. Hence
\begin{eqnarray*}
\sup\{\frac{|q(y)|}{\omega^*(y)}: y\in H\}&\geq&\sup\{\frac{|q(\tilde{x})|}{\omega^*(\tilde{x})}: x\in G\}\\
&=&\sup\{\frac{|Q( x)|}{\omega^*(\tilde{x})}: x\in G\}=\sup\{\frac{|Q(x)|}{(\tilde{\omega})^*(x)}: x\in G\}.
\end{eqnarray*}
Therefore, $L^1(H, \omega|_H)$ is weakly amenable.$\hfill\square$

%begin{example}{\rm Let $G={\Bbb Z}$, $H=2\;{\Bbb Z}$ and $\omega(n)=\omega_{1/3}(n)$ for all $n\in{\Bbb Z}$. Then
%$G=H\cup (1+H)$. So $A=\{0, 1\}$. Note that
%\begin{eqnarray*}
%\tilde{n}=\left\{
%\begin{array}{rl}
%n & n\;\; \hbox{is even}\\
%-n    & n\;\; \hbox{is odd}
%\end{array}\right.
%\end{eqnarray*}
%and hence
%\begin{eqnarray*}
%\tilde{\omega}(n)=\left\{
%\begin{array}{rl}
%\omega_{1/3}(n) & n\;\; \hbox{is even}\\
%\omega_{1/3}(-n)    & n\;\; \hbox{is odd}
%\end{array}\right.
%\end{eqnarray*}}
%\end{example}

%\begin{proposition} Let $N$ be a normal subgroup of a locally compact group $G$ and $\omega$ be a weight function on $G$.
%\emph{(i)} If $L^1(G, \omega)$ is weakly amenable, then $L^1G/N, \omega^\prime)$ is weakly amenable.
 %If $L^1(G\times G, \omega^\otimes)$ is amenable, then
%$L^1(N\times G/N, \omega|_H\otimes\hat{\omega})$ is weakly amenable.
%\end{proposition}
%{\it Proof.}
% We only note that  $$\omega|_H\otimes\hat{\omega}(h, xH)\leq\omega(h)\omega(x)$$
%for all $x\in G$ and $h\in H$.$\hfill\square$\\

Now, we investigate weak amenability of tensor product of Beurling algebras.

\begin{theorem}\label{ppp} Let $G_1$ and $G_2$ be locally compact groups and let $\omega_1$ and $\omega$ be weight functions on $G_1$ and $G_1\times G_2$, respectively. If $L^1(G_1\times G_2, \omega)$ is weakly amenable and $\sup\{\frac{\omega_1(x)}{\omega(x, y)}: x, y\in G\}<\infty$, then $L^1(G_1, \omega_1)$ is weakly amenable.
\end{theorem}
{\it Proof.} Let $p$ be a non-inner quasi-additive function in  $L^\infty(G_1, 1/\omega_1)$. Define the function  $\tilde{p}((x_1, y_1), (x_2, y_2))=p(x_1, x_2)$ for all $x_1, y_1\in G_1$ and $x_2, y_2\in G_2$. Suppose that $\tilde{p}\in I(G_1\times G_2, \omega_1\otimes\omega_2)$. Then there exists $h\in L^\infty(G_1\times G_2, \omega_1\times\omega_2)$ such that
$$
\tilde{p}(x_1,x_2),(y_1,y_2))=h((x_1, y_1)(x_2, y_2))-h((x_2, y_2)(x_1, y_1))
$$
for all $x_1, y_1\in G_1$ and $x_2, y_2\in G_2$. Define the complex-valued function $k$ on $G_1$  by $k(x)=h(x, e_2)$ for all $x\in G_1$.
Note that
$$
\frac{k(x_1)}{\omega_1(x_1)}=\frac{h(x_1, e_2)}{\omega_1\otimes\omega_2(x_1, e_2)}.
$$
This shows that $k\in L^\infty(G_1, 1/\omega_1)$. For every $x_1, y_1\in G_1$ and $x_2, y_2\in G_2$, we have
\begin{eqnarray*}
k(x_1y_1)-k(y_1x_1)&=&h((x_1, e_2)(y_1, e_2))-h((y_1, e_2)(x_1, e_2))\\
&=&\tilde{p}((x_1, e_2), (y_1, e_2))\\
&=&p(x_1, y_1).
\end{eqnarray*}
This contradiction shows that $\tilde{p}$ is non-inner. If $x_1, y_1\in G_1$ and $x_2, y_2\in G_2$,then
$$
\frac{|\tilde{p}(x_1,x_2),(y_1,y_2))|}{\omega(x_1,x_2)\omega(y_1,y_2)}=
\frac{|p(x_1,y_1)|}{\omega_1(x_1)\omega_1(y_1)}\frac{\omega_1( x_1)}{\omega( x_1,x_2)}
\frac{\omega_1(y_1)}{\omega(y_1,y_2)}.
$$
This implies that if $L^1(G_1, \omega_1)$ is not weakly amenable, then $L^1(G_1\times G_2, \omega)$ is not weakly amenable.$\hfill\square$\\

For a locally compact Abelian group $G$, Zhang \cite{z} proved that if $L^1(G_1, \omega_1)\hat{\otimes} L^1(G_2, \omega_2)$ is weakly amenable, then $L^1(G_1, \omega_1)$ and $L^1(G_2, \omega_2)$ are weakly amenable. We show that this result is true for any locally compact group.

\begin{corollary}\label{ten} Let $\omega_i$ be a weight function on a locally compact infinite group $G_i$, for $i=1, 2$. If $L^1(G_1, \omega_1)\hat{\otimes} L^1(G_2, \omega_2)$ is weakly amenable, then $L^1(G_1, \omega_1)$ and $L^1(G_2, \omega_2)$ are weakly amenable.
\end{corollary}
{\it Proof.} Let $L^1(G_1, \omega_1)\hat{\otimes} L^1(G_2, \omega_2)$ be weakly amenable. Then $L^1(G_1\times G_2, \omega_1\otimes\omega_2)$ is weakly amenable. Since for every $x, y\in G$
$$
\frac{\omega_1(x)}{\omega_1\otimes\omega_2(x, y)}=\frac{1}{\omega_2(y)}\leq 1,
$$
it follows from Theorems \ref{ppp} that $L^1(G_1, \omega_1)$ is weakly amenable. The other case is similar.$\hfill\square$

%\begin{example}{\rm Let $G_1$ and $G_2$ be locally compact groups and let $\omega$ is a weight function on $G_1\times G_2$. Define the weight functions $\omega_1$ and $\omega_2$ on $G_1$ and $G_2$, respectively, by
%$$
%\omega_1(x_1):=\omega(x_1, e_2)\quad\hbox{and}\quad\omega_2(x_2):=\omega(e_1, x_2)
%$$
%for all $x_1\in G_1$ and $x_2\in G_2$. It is clear that $\omega\leq\omega_1\otimes\omega_2$. So if $L^1(G_1\times G_2, \omega_1\otimes\omega_2)$ is weakly amenable, then $L^1(G_1\times G_2, \omega)$ is weakly amenable.}
%\end{example}

\vspace{2mm}

 {\footnotesize
\noindent {\bf Mohammad Javad Mehdipour}\\
Department of Mathematics,\\ Shiraz University of Technology,\\
Shiraz
71555-313, Iran\\ e-mail: mehdipour@.ac.ir\\
{\bf Ali Rejali}\\
Department of Pure Mathematics,\\ Faculty of Mathematics and Statistics,\\ University of Isfahan,\\
Isfahan
81746-73441, Iran\\ e-mail: rejali@sci.ui.ac.ir\\

\footnotesize

\begin{thebibliography}{99}

\bibitem{bcd} W. Bade, P. C. Curtis and H. G. Dales, Amenability and weak amenability for Beurling and Lipschitz algebras, Proc. London Math. Soc., (3) 55 (1987) 359--377.

\bibitem{bd} F. F. Bonsall and J. Duncan, Complete Normed Algebras, Ergebnisse der Mathematik und ihrer Grenzgebiete, Band 80, Springer-Verlag, Berlin/ Heidelberg/New York, 1973.

\bibitem{b} C. R. Borwick, Johnson-Hochschild Cohomology of Weighted Group Algebras and Augmentation Ideals, Ph. D. thesis, University of Newcastle upon Tyne, 2003.

\bibitem{d} H. G. Dales, Banach algebras and Automatic Continuity, Clarendon Press, Oxford, 2000.

\bibitem{dl} H. G. Dales and A. T. Lau, The second duals of Beurling algebras, Mem. Amer. Math. Soc., 177 (836) (2005).

\bibitem{dg} M. Despic and F. Ghahramani, Weak amenability of group algebras of locally compact groups, Canad. Math. Bull., 37 (2) (1994) 165--167.

\bibitem{gro} N. Gronback, A charactrization of weakly amenable Banach algebras, Studia. Math., 94 (1989) 149--162.

\bibitem{hr} E. Hewitt and K. Ross,
Abstract Harmonic Analysis I, Springer-Verlag, New York, 1970.

\bibitem{j} B. E. Johnson, Weak amenability of group algebras, Bull. London Math. Soc., 23 (3) (1991) 281--284.

\bibitem{mr1} S. Maghsoudi and A. Rejali, Unbounded weighted Radon measures and dual of certain function spaces with strict topology, Bull. Malays. Math. Sci. Soc., 36 (1) (2013) 211--219.

\bibitem{mr2} S. Maghsoudi and A. Rejali, On the dual of certain locally convex function spaces, Bull. Iranian Math. Soc., 41 (4) (2015) 1003--1017.

\bibitem{mr} M. J. Mehdipour and A. Rejali, Regularity and amenability of weighted Banach algebras and their second dual on locally compact groups,  arXiv:2112.13286v1.

\bibitem{p} A. Pourabbas, Weak amenability of weighted group algebras, Atti Sem. Mat. Fis. Univ. Moderna, 48 (2) (2000) 299--316.

\bibitem{r0} H. Reiter and J. D. Stegeman, Classical Harmonic Analysis and Locally Compact Groups, London Math. Society Monographs, 22, Clarendon Press, Oxford, 2000.

\bibitem{r111} A. Rejali, The analogue of weighted group algebra for semitopological semigroups, J. Sci. Islam. Repub. Iran, 6 (2) (1995) 113--120.

\bibitem{rv1} A. Rejali and H. R. Vishki, Weighted convolution measure algebras characterized by convolution algebras, J. Sci. Islam. Repub. Iran, 19 (2) (2008) 169--173.

\bibitem{s1} V. Shepelska, Weak amenability of weighted group algebras on some discrete groups, Studia Math., 230 (3) (2015) 189--214.

\bibitem{sto} R. Stokke, On Beurling measure algebras, arXiv:2107.14694v1.

\bibitem{z} Y. Zhang, Weak amenability of commutative Beurling algebras, Procc. Amer. Math. Soc., 142 (5) (2014) 1649--1661.
\end{thebibliography}
\end{document}